\documentclass[a4paper]{article}

\usepackage{hyperref}

\usepackage{amsmath,amssymb,amsthm}
\usepackage{cases}
\numberwithin{equation}{section}

\usepackage{shortvrb, makeidx}

\makeindex

\title{Formulation of finite-time singularity for free-surface Euler equations}
\author{Yi Zhou\thanks{  School of Mathematical Sciences, Fudan
University, Shanghai 200433, P. R. China ({\tt Email:
yizhou@fudan.ac.cn})
  }}
\date{}
\begin{document}
\maketitle
\begin{abstract}
We give an extremely short proof that the free-surface incompressible, irrotational Euler equations with regular initial condition can form a finite time singularity in 2D or 3D. Thus, we provide a simple view of the problem studied by Castro et al \cite{Castro1}, \cite{Castro2}, \cite{Castro3}, \cite{Castro4},
\cite{Castro5} and Coutand \& Shkoller \cite{Coutand3}.
\end{abstract}
\section{\textbf{Introduction}}

\indent The aim of this paper is to give an extremely short proof that the free-surface incompressible, irrotational Euler equations with regular  initial conditions can form a finite time singularity in 2D or 3D. Before going further, we make some historical remarks on the problem. For the irrotational case of the water waves problem, the local wellposedness in Sobolev space was established by Wu \cite{Wu1} in 2D and Wu \cite{Wu2} in 3D. Earlier works dealt with small data or linearized equations, see Nalimov \cite{Nalimov}, Yoshihara \cite{Yosihara}, Craig \cite{Craig}, Beale et al \cite{Beale}. Further studies of local wellposedness was done by Ambrose \&Masmoudi \cite{Ambrose1} and by Lannes \cite{Lannes}. For small initial data, Wu \cite{Wu3} established almost global existence for the infinite-depth problem in 2D. Independently, Wu \cite{Wu4} and Germain et al \cite{Germain1}, \cite{Germain3} established global existence for the infinite-depth problem in 3D. For more bibliographic notes, we refer to \cite{Christodoulou}, \cite{Coutand1}, \cite{Coutand2}, \cite{Germain2}, \cite{Lindblad}, \cite{Shatah}, \cite{Zhang}.\\
\indent The problem of establishing a finite-time singularity for the interface has recently been explored for the 2D water waves equations by Castro et al \cite{Castro3}, \cite{Castro4}, where it was shown that a smooth initial curve exhibits a finite-time singularity via self-intersection at a point. More recently, Coutand and Shkoller \cite{Coutand3} established a similar result in 3D for the general incompressible case, and their method can also be used in 2D. Both approaches are to start from a curve which self-intersecting at a point and solve backward to get a regular solution at an early time. In this paper,  we start from given initial conditions and show by contradiction that a finite-time singularity must form if the initial data satisfies some positivity condition. Thus, we provide an extremely simple view of the problem. However, the price is to fix the endpoint as will be see from below. Our method works both in 2D and 3 D, with or without gravity, for simplicity, we just give a proof in 2D without gravity. As our proof is extremely simple, it is not hard to figure out a same proof in 3D or with gravity. \\
\indent We consider the 2D free-surface incompressible, irrotational Euler equations in $\Omega(t)$, where $\Omega(t)$ is bounded by $x_1=0,1,$ and $x_2=0$ and a free surface $\Gamma(t)$. We assume that $\Gamma(t)$ is parametered as
\begin{equation}
\Gamma (t)=\{\big(x_1(t,\alpha),x_2(t,\alpha)\big),0\leq\alpha\leq1\},
\end{equation}
with \begin{equation}
x_2(t,0)=x_2(t,1)=1,
\end{equation}
\begin{equation}
0\leq x_1(t,\alpha)\leq 1, 0\leq t<\infty, 0\leq \alpha \leq1.
\end{equation}
Our problem is modeled by the following
\begin{equation}
\left \{
\begin{aligned}
&u^1_t+u^1\partial_1u^1+u^2\partial_2u^1+\partial_1p=0,\\
&u^2_t+u^1\partial_1u^2+u^2\partial_2u^1+\partial_2p=0,~~ on~ \Omega(t)\\
\end{aligned} \right.
\end{equation}
\begin{equation}
\left \{
\begin{aligned}
&\partial_1u^1+\partial_2u^2=0,\\
&\partial_2u^1-\partial_1u^2=0,~~ on~ \Omega(t)\\
\end{aligned} \right.
\end{equation}
\begin{equation}
\left \{
\begin{aligned}
&u^1=0,~~ ~ x_1=0,1\\
&u^2=0,~~~x_2=0,~~ on~ \Omega(t)\\
\end{aligned} \right.
\end{equation}
\begin{equation}
u^2(t,0,1)=u^2(t,1,1)=0
\end{equation}
\begin{equation}
\left \{
\begin{aligned}
&\frac{dx_1(t,\alpha)}{dt}=u^1(t,x_1(t,\alpha),x_2(t,\alpha)),~~0\leq\alpha\leq1\\
&\frac{dx_2(t,\alpha)}{dt}=u^2(t,x_1(t,\alpha),x_2(t,\alpha)),\\
\end{aligned} \right.
\end{equation}
\begin{equation}
t=0: u^1=u^1_0, u^2=u^2_0 ~~on~~\Omega(0)
\end{equation}
where $u^1_0, u^2_0 $ satisfies (1.5) and (1.2), (1.3), (1.6), (1.7) at $t=0$.\\
Our main result can be summarized as follows:\\
{\bf Theorem 1.1.} Consider the free-surface Euler equation (1.1)-(1.9), let $\Gamma(t)$ and $u^1, u^2$ be a $C^1$ solution with $C^1$ initial data $\Gamma(0)$ and $(u^1_0, u^2_0)$, then $\Gamma(t)$ must develop a finite time singularity provided that
$$A=\int_{\Omega(0)}u^1_0 x_1dx+\int_0^1x_2u^2_0(x_2,1)dx_2>0.$$
\section{\textbf{Proof of main result}}
As observed by Lindblad \cite{Christodoulou}
\begin{equation}
\begin{aligned}
-\Delta p &=(\partial_1u^1)^2+(\partial_2u^2)^2+2\partial_1u^2\partial_2u^1\\
&=(\partial_1u^1)^2+(\partial_2u^2)^2+2(\partial_2u^1)^2\\
&\geq 0,
\end{aligned}
\end{equation}
where $\triangle=\partial^2_1+\partial^2_2$,
$$p=0~~~~on ~\Gamma(t),$$
it is not difficult to see from boundary condition that
\begin{equation}
\frac{\partial p}{\partial n}=0, ~x_1=0,1~or~x_2=0,
\end{equation}
where $n$ is outward norm.\\
Therefore, by maximum principle, we get
\begin{equation}
p\geq0~ on ~\Omega(t).
\end{equation}
Let
\begin{equation}
\frac{D}{Dt}=\partial_t+u^1\partial_1+u^2\partial_2
\end{equation}
be the material derivative, then
\begin{equation}\begin{aligned}
\frac{D(u^1x_1)}{Dt}&=\frac{Du_1}{Dt}x_1+u_1\frac{Dx_1}{Dt}\\&
=-\partial_1p~x_1+u_1^2.
\end{aligned}
\end{equation}
Therefore
\begin{equation}\begin{aligned}
\frac{d}{dt}\int_{\Omega(t)}u^1x_1dx&=\int_{\Omega(t)}(u^1)^2dx-\int_{\Omega(t)}\partial_1p~x_1dx\\
&=\int_{\Omega(t)}(u^1)^2dx+\int_{\Omega(t)}p dx-\int_0^1p(t,1,x_2)dx_2
\end{aligned}
\end{equation}
On $x_1=1, u^1=0$, so we get
$$\partial_t u^2+u^2\partial_2 u^2+ \partial_2 p=0,$$
multiply by $x_2$ and make an integration by parts, we get
\begin{equation}
\frac{d}{dt}\int_0^1u^2(t,x_2,1)x_2dx_2=\frac{1}{2}\int(u^2(t,x_1,1))^2dx_2+\int_0^1p(t,1,x_2)dx_2.
\end{equation}
Adding (2.6) and (2.7) together and noting (2.3), we get
\begin{equation}\begin{aligned}
\frac{d}{dt}&\int_{\Omega(t)}u^1x_1dx+\int_0^1u^2(t,x_2,1)x_2dx_2\\
&=\int_{\Omega(t)}(u^1)^2dx+\frac{1}{2}\int_0^1(u^2(t,x_2,1))^2dx_2+\int_{\Omega(t)}pdx\\
&\geq\int_{\Omega(t)}(u^1)^2dx+\frac{1}{2}\int_0^1(u^2(t,x_2,1))^2dx_2.
\end{aligned}
\end{equation}
By Schwatz inequality and noting (1.5) and (1.3), we get
\begin{equation*}\begin{aligned}
\left(\int_{\Omega(t)}u^1x_1dx \right)^2&\leq\left(\int_{\Omega(t)}(u^1)^2dx\right)\left(\int_{\Omega(t)}x_1^2dx\right)\\
&\leq\left(\int_{\Omega(t)}(u^1)^2dx\right)\int_{\Omega(t)}dx\\
&=\int_{\Omega(t)}(u^1)^2dx\int_{\Omega(0)}dx.\\
\end{aligned}
\end{equation*}
Similarly
\begin{equation*}\begin{aligned}
\big(\int_0^1x_2&u^2(t,x_2,1)dx_2\big)^2\\
&\leq \left(\int_0^1x_2^2dx_2\right)\int_0^1\left(u^2(t,x_2,1)\right)^2dx_2\\
&=\frac{1}{3}\int_0^1\big(u^2(t,x_2,1)\big)^2dx_2.
\end{aligned}
\end{equation*}
Let
\begin{equation}
c_1=max\big(2\int_{\Omega(0)}dx,~\frac{4}{3}\big)
\end{equation}
and
\begin{equation}
L(t)=\int_{\Omega(t)}u^1x_1dx+\int_0^1x_2u^2(t,x_2,1)dx_2
\end{equation}
then it is not difficult to see
\begin{equation}\label{eq3}
\left \{
\begin{aligned}
&L^\prime(t)\geq\frac{L(t)^2}{c_1}\\
&L(0)=A>0
\end{aligned} \right.
\end{equation}
Thus, $L(t)$ becomes infinite in finite time.
\section*{Acknowledgments}
This research is supported by Key Laboratory of Mathematics for Nonlinear Sciences (Fudan University), Ministry of Education of China, P. R. China, Shanghai Key Laboratory for Contemporary Applied Mathematics, School of Mathematical Sciences, Fudan University,  NSFC(grants No. 11031001 and 11121101), and 111 project.


\begin{thebibliography}{4}
\addcontentsline{toc}{section}{References}
\bibitem{Alvarez-Samaniego}
B. Alvarez-Samaniego and D. lannes, Large time existence for 3D water-waves and asymptotics. Invent. Math., 171, (2008), 485-541.
\bibitem{Ambrose1}
D. M. Ambrose and N. Masmoudi, The zero surface tension limit of two-dimensional water waves, Comm. Pure. Appl. Math., 58(2005), 1287-1315.
\bibitem{Ambrose2}
D. M. Ambrose and N. Masmoudi, The zero surface tension limit of three-dimensional water waves, Indiana Univ. Math. J., 58(2009), 479-521.
\bibitem{Beale}
J. T. Beale, T. Hou and J. Lowengrub, Growth rates for the linearized motion of fluid interfaces away from equilibrium,Comm. Pure Appl. Math.,46(1993), 1269-1301.
\bibitem{Castro1}
A. Castro, D. C$\acute{o}$rdoba, C. Fefferman, F. Gancedo and M. L$\acute{o}$pez-Fern$\acute{a}$ndez, Rayleigh-Taylor breakdown for the Muskat problem with applications to water waves, Ann. of Math., 2,(2011), to appear.
\bibitem{Castro2}
A. Castro, D. C$\acute{o}$rdoba, C. Fefferman, F. Gancedo and M. L$\acute{o}$pez-Fern$\acute{a}$ndez, Turning waves and breakdown for incompressible flows, Proceedings of the National Academy of Sciences, 108(2011), 4754-4759.
\bibitem{Castro3}
A. Castro, D. C$\acute{o}$rdoba, C. Fefferman, F. Gancedo and M. G$\acute{o}$mez-Serrano, Splash singularity for water waves, (2011), arxiv: 1106.2120v2.
\bibitem{Castro4}
A. Castro, D. C$\acute{o}$rdoba, C. Fefferman, F. Gancedo and M. G$\acute{o}$mez-Serrano,  Finite time singularities for the free boundary incompressible Euler equaitons, (2011), arxiv: 1112.2170v1.
\bibitem{Castro5}
A. Castro, D. C$\acute{o}$rdoba, C. Fefferman, F. Gancedov and M. G$\acute{o}$mez-Serrano, Finite time singularities for water waves with surface tension, preprint 2012.
\bibitem{Craig}
W. Craig, An existence theory for water waves and the Boussinesq and Korteweg-de Vries scaling limits. Comm. Partial Differential Equations, 10(1985),no.8, 787-1003.
\bibitem{Christodoulou}
D. Christodoulou and H. Lindblad, On the motion of the free surface of a liquid, Comm. Pure Appl. Math., 53(2000), 1536-1602.
\bibitem{Coutand1}
D. Coutand and S. Shkoller, Well-posedness of the free-surface incompressible Euler equations with or without surface tension, J. Amer. Math. Soc.,20(2007), 829-930.
\bibitem{Coutand2}
D. Coutand and S. Shkoller, A simple proof of well-posedness for the free-surface incompressible Euler equations, Discrete Contin. Dyn. Syst. Ser. S, 3 (2010), 429-449.
\bibitem{Coutand3}
D. Coutand and S. Shkoller, On the finite-time splash singularity for the 3-D free-surface Euler equations, preprint, 2012.
\bibitem{Germain1}
P. Germain, N. Masmoudi and J. Shatah, Global solutions for the gravity water equation in dimension 3, C. R. Math. Acad. Sci. Paris, 347(2009), 897-902.
\bibitem{Germain2}
P. Germain, N. Masmoudi and J. Shatah, Global existence for capillary water waves, preprint 2012.
\bibitem{Germain3}
P. Germain, N. Masmoudi and J. Shatah, Global solutions for the gravity water waves equation in dimension 3, Ann. of Math, to appear
\bibitem{Lannes}
D. Lannes, Well-posedness of the water-waves equations, J. Amer. Math. Soc., 18(2005), 605-654.
\bibitem{Lindblad}
H. Lindblad, Well-posedness for the motion of an incompressible liquid with free surface boundary, Ann. of Math., 162(2005), 109-194.
\bibitem{Nalimov}
V.I. Nalimov, The Cauchy-Poisson Problem (in Russion), Dynamika Splosh. Sredy, 18(1974), 104-210.
\bibitem{Shatah}
J. Shatah and C. Zeng, Geometry and a priori estimate for free boundary problems of the Euler equation, Comm. Pure Appl. Math., 61(2008), 698-744.
\bibitem{Wu1}
S. Wu, Well-posedness in Soblev spaces of the full water wave problem in 2-D, Invent. Math., 130(1997),39-72.
\bibitem{Wu2}
S. Wu, Well-posedness in Soblev spaces of the full water wave problem in 3-D, J. Amer. Math. Soc., 12 (1999), 445-495.
\bibitem{Wu3}
S. Wu, Almost global wellposedness of the 2-D full water wave problem, Invent. Math., 177(2009),39-72.
\bibitem{Wu4}
S. Wu, Global wellposedness of the 3-D full water wave problem, Invent. Math., 184(2011),125-220.
\bibitem{Yosihara}
H. Yosihara, Gravity Waves on the Free Surface of an Incompressible Perfect Fluid, Publ. RIMS Kyoto Uni., 18(1982), 49-96.
\bibitem{Zhang}
P. Zhang and Z. Zhang, On the free boundary problem of three-dimensional incompressible Euler eqautions, Comm. Pure Appl. Math., 61(2008), 877-940.
\end{thebibliography}
\end{document}